%%%%%%%%%%%%%%%%%%%%%%%%%%%%%%%%%%%%%%%%%%%%%%%%%%%%%%%%%%%%%%%%%%%%%%%%%%%
%% de la Pena, Victor H.; Montgomery-Smith, S. J.
%% 
%% Bounds on the tail probability of U-statistics and quadratic forms
%% 
%% publ:  Bull. Amer. Math. Soc. (N.S.) 31(1994) no. 2
%% pp:    223-227
%% type:  Research Announcement        markup: amstex    file size: 18K
%% contact: vp@wald.stat.columbia.edu
%% 
%% copyright: American Math. Society copyright; see end of article
%% 
%% Include files necessary for this article: bull-ppt.tex
%% 
%%%%%%%%%%%%%%%%%%%%%%%%%%%%%%%%%%%%%%%%%%%%%%%%%%%%%%%%%%%%%%%%%%%%%%%%%%%
% Date: 9-MAR-1994
%   : 0   \newpage: 0   \displaybreak: 0
%   \eject: 0   \bye: 0   \break: 0   \allowbreak: 0
%   \allowdisplaybreak: 0   \allowdisplaybreaks: 0
%   \allowlinebreak: 0   \allowmathbreak: 0
%   \smallpagebreak: 0   \medpagebreak: 0   \bigpagebreak: 0
%   \smallbreak: 0   \medbreak: 0   \bigbreak: 0   
%\goodbreak: 0
%   : 0   : 0   \newline: 0
%   \magnification: 1   \mag: 0
%   \baselineskip: 0   \normalbaselineskip: 0
%   \hsize: 0   \vsize: 0   \pagewidth: 0   \pageheight: 0
%   \hoffset: 0   \voffset: 0   \hcorrection: 0   
%\vcorrection: 0
%   \parindent: 0   \parskip: 0
%   \vfil: 0   \vfill: 0   \vskip: 0
%   \smallskip: 0   \medskip: 0   \bigskip: 1
%   \sl: 0   \def: 0   \let: 0   \redefine: 0   
%\predefine: 0
%   \tolerance: 0   \pretolerance: 0
%   \font: 0   \end: 0   \\noindent: 3
%   ASCII 13 (Control-M Carriage return): 0
%   ASCII 10 (Control-J Linefeed): 0
%   ASCII 12 (Control-L Formfeed): 0
%   ASCII 0 (Control-@): 0
% Special characters: 0
%
%
\input amstex
\documentstyle{amsppt}
\input bull-ppt
\keyedby{pena-1/amh}   % (xxxxxx is the paper no.; ### is 
%the user's initials)
%\documentstyle{amsppt}
%\magnification=1200

\topmatter
\cvol{31}
\cvolyear{1994}
\cmonth{October}
\cyear{1994}
\cvolno{2}
\cpgs{223-227}
\title Bounds on the tail probability \\
of  U-statistics and quadratic forms \endtitle
\author Victor H. de la Pe\~na and S.~J.~Montgomery-Smith 
\endauthor
\shortauthor{V. H. de la Pe\~na and S.~J.~Montgomery-Smith}
\address Department of Statistics, Columbia University, 
New York, New York 10027\endaddress %\endaddress
\ml %\email 
vp\@wald.stat.columbia.edu\endml %\endemail
\address Department of Mathematics, University of 
Missouri, Columbia, Missouri
65211\endaddress %\endaddress
%\email 
\date August 3, 1993, and, in revised form, August 24, 
1993\enddate
\ml stephen\@mont.cs.missouri.edu\endml %email
\thanks Both authors were supported by NSF grants. The 
second author was
supported by the Research Council of the University of 
Missouri\endthanks
 % \endsubjclass
%\endkeywords
\endtopmatter

\document
It is very common for expressions of the form:
$$ \sum_{1\le i_1 \ne \cdots \ne i_k \le n}
   f_{i_1 \cdots i_k}(X_{i_1},\dots,X_{i_k}) $$
to appear in probability theory.  Here $\{X_i\}$\ is a 
sequence of independent
random variables taking values in a measurable space $(S, 
{\Cal S})$,
and $\{f_{i_1\cdots i_k}\}$ is a sequence of measurable 
functions
from  $S^k$ into a 
Banach Space $(B,
\| \cdot \|)$.  Special cases of this type of random 
variable appear, 
for example, in statistics in
the form of U-statistics and quadratic forms.
Throughout we will refer to them as generalized 
U-statistics.

There is great interest in decoupling such quantities, 
that is, in replacing
the above quantity by the expression
$$ \sum_{1\le i_1 \ne \cdots \ne i_k \le n}
   f_{i_1 \cdots i_k}(X_{i_1}^{(1)},\dots,X_{i_k}^{(k)}),$$
where $\{X_i^{(1)}\}$, $\{X_i^{(2)}\},\dots$,
$\{X_i^{(k)}\}$ are $k$\ independent
copies of $\{X_i\}$.  

Decoupling inequalities allows one to compare expressions 
of the first kind
with expressions of the second kind.  
Such results permit the almost-direct transfer of results 
for sums of
independent random variables to the case of generalized 
U-statistics. 
The reason for this is that,
conditionally on $\{X_i^{(2)}\}\!,...\!,\{X_i^{(k)}\}$\<,
the second sum above is a sum of 
independent random variables. 
It is important to remark that such results have 
led to the development of several optimal 
results in the functional theory of U-statistics (cf. [1]
and [7]) and various other areas,
including the study of the invertibility of large matrices
(cf. [2]),  
stochastic integration (cf. [10]),
and the study of integral operators on Lebesgue-Bochner 
spaces 
(cf. a result of T. R. McConnell and D. Burkholder found in
[3]). Aside from those directly cited in this paper, other
important contributors to the area of
decoupling inequalities include
A. de Acosta,  
P. Hitczenko, 
J. Jacod, 
A. Jakubowski, 
O. Kallenberg, 
M. Klass,
W. Krakowiak,
G. Pisier,
J. Rosinski, and
J. Szulga.
Due to space restrictions, we refer the reader to 
[10] for a more complete account.

In this paper, we announce a result which allows one to 
compare
the tail probabilities of the above quantities.
In particular, this inequality 
represents the definitive generalization 
of the decoupling inequalities for multilinear forms of 
McConnell
and Taqqu [11] and the more general decoupling 
inequalities for
expectations of convex functions of U-statistics 
introduced in 
[4]. 
\proclaim{Theorem 1} There is a constant $C_k>0$, 
depending only on $k$, 
such that
for all $n\ge k$, 
$$\split
& P\left(\left\|\sum_{1\le i_1 \ne \dots \ne i_k \le n}
   f_{i_1 \dots i_k}(X_{i_1},\dots,X_{i_k}) \right\|\ge 
t\right) \\
&\qquad    \leq C_k P\left(C_k
\left\|\sum_{1\le i_1 \ne \dots \ne i_k \le n}
   f_{i_1 \dots i_k}(X_{i_1}^{(1)},\dots,X_{i_k}^{(k)}) 
\right\|\ge t\right),
 \quad\text{ for all $t>0$.}\endsplit\tag1$$
Moreover, the reverse inequality holds if 
the functions satisfy the condition
$$ f_{i_1 \dots i_k}(X_{i_1},\dots,X_{i_k}) =
   f_{i_{\pi(1)} \dots 
i_{\pi(k)}}(X_{i_{\pi(1)}},\dots,X_{i_{\pi(k)}}) $$
for all permutations $\pi$\ of $\{1,\dots,k\}$. 
That is, 
$$\split
& P\left(\left\|\sum_{1\le i_1 \ne \dots \ne i_k \le n}
   f_{i_1 \dots i_k}(X_{i_1},\dots,X_{i_k}) \right\|\ge 
t\right) \\
&\qquad \geq 
  {1\over C_k} P\left({1\over C_k}
\left\|\sum_{1\le i_1 \ne \dots \ne i_k \le n}
   f_{i_1 \dots i_k}(X_{i_1}^{(1)},\dots,X_{i_k}^{(k)}) 
\right\|\ge t\right),
\quad \text{ for all $t>0$.}\endsplit\tag2$$
Note that
the expression $  i_1 \ne \dots \ne i_k $\ means that $i_r 
\ne i_s$\ for all
$1\le r \ne  s \le k$.
\endproclaim

An example illustrating this result can be found in the 
study of random
graphs (see also [5]).  
Given a sequence of independent random points $\{X_i\}$\
in $R^N$, we might consider a measure of clustering
$$ D_1 = \sum_{1\le i \ne j \le n} d(X_i,X_j) ,$$
where $d(x,y)$\ denotes the distance between $x$\ and $y$.
The above result allows us to compare $D_1$, which 
measures the distance
``within'' the graph formed by the random cluster of points 
$\{X_i\}$, to a quantity $D_2$, which is a  measure of the 
distance ``between''
the two independent clusters $\{X_i\}$ and $\{\tilde X_i 
\}$,
$$D_2 = \sum_{1\le i \ne j \le n} d(\tilde X_i, X_j), $$
%\noindent 
where $\{\tilde X_i\}$ is an independent copy of $\{X_i\}$. 
Then we have for all $t>0$ that
$$ C_2^{-1} P(|D_1| \ge C_2 t) \le P(|D_2| \ge t)
   \le C_2 P(|D_1| \ge C_2^{-1} t) .$$

Other examples where U-statistics are used in graph theory 
may be found
in [8]. 

%\bigskip
We will prove the theorem in the special case that $k=2$. 
For ease of notation, let us suppose that $ \tilde X_i = 
X_i^{(1)}$
and denote $ X_i = X_i^{(2)}$.
The proof
of the more general result will appear elsewhere.
We will use a sequence of lemmas. 
Following [6], 
our point of departure is equation (4), which
provides a partial decoupling result and focuses attention 
on a
polarized version of the U-statistic kernel as the key 
element in 
the development of a solution of the problem at hand. 
Let 
$$T_n =
\sum_{1\le i\ne j\le n} \{f_{ij}(X_i, X_j)+
f_{ij}(X_i, \tilde X_j)+   
f_{ij}(\tilde X_i, X_j)
+ f_{ij}(\tilde X_i,\tilde X_j)\}; \leqno (3)$$ then by 
using the
triangle inequality, one obtains that
$$\split
&P\left(\left\|\sum_{1\le i\ne j\le n} f_{ij}(X_i, X_j) 
+ f_{ij}(\tilde X_i,\tilde X_j)\right\|\ge t \right)\\
&\qquad \le 
P(\| T_n \| \ge { t \over 3} ) +
2P\left(
\left\|\sum_{1\le i\ne j\le n} f_{ij}(X_i, \tilde X_j)
\right\| \ge { t \over 3} \right)\.  \endsplit\tag4 $$
This observation reduces the proof of (1) to the problem  
of 
obtaining the bounds 
$$\aligned &P\left(\left\|\sum_{1\le i\ne j\le n} 
f_{ij}(X_i, X_j)
\right\|\ge t\right)\\&\qquad \le
cP\left(c\left\|\sum_{1\le i\ne j\le n} f_{ij}(X_i, X_j) +
f_{ij}(\tilde X_i,\tilde X_j)\right\|\ge 
t\right)\!,\endaligned\leqno (5)$$
and
$$  P(\|T_n
\| \ge { t } ) 
 \le  
cP\lf(c\lf\|\sum_{1\le i\ne j\le n} f_{ij}
(\tilde X_i,X_j)\rt\|  \ge  t\rt).\leqno (6)$$ 
We obtain (5) by means of Lemma 1 
(possibly of independent interest).
The proof of (6) is somewhat involved. In obtaining it, we 
used 
(conditionally) an extension of the Paley-Zygmund 
inequality  found in
[10] in combination with a symmetrization identity 
similar to the one introduced in [12].
\proclaim{Lemma 1} Let X, Y be two i.i.d. random variables.
Then
$$
P(\|X\|\ge t)  \le 
3P(\|X+Y\|\ge {2t\over 3}). \leqno (7) $$
\endproclaim
\demo{Proof}
Let $X$, $Y$, and $Z$ be i.i.d. random
variables. Then
$$\eqalign{P(\|X \|& \ge t)  \cr
&=
P(\|(X+Y) + (X+Z) - (Y+Z) \ge 2t) \cr
& \le
P(\|X+Y\| \ge {2t\over 3})+P(\|X+Z\| \ge {2t\over 3})+
P(\|Y+Z\| \ge {2t\over 3}) \cr
& =
3P(\|X+Y\| \ge {2t\over 3}). \cr }$$
\enddemo

It is to be remarked that the
very desirable ``Universal Symmetrization 
Lemma'', 
$P(\|X\|\ge t)  \le 
c P(c \|X-Y\|\ge t)$, is not true.  This makes the above 
result all the more
surprising.

The following is an observation found in Section 6.2 of
[10] that will be used in combination with 
Lemma 2 to prove Theorem 1.
\proclaim{Proposition 1}
Let $Y$ be any mean-zero random variable with values in a
Banach space $(B, \|\cdot \|)$. Then, for all $a\in B$,
$P(\|a+Y\| \ge \|a\|) \ge {\kappa \over 4}$, where,
$\kappa = \inf_{x'\in B'}
 {(E|x'(Y)|)^2 \over E(x'(Y))^2}$.
\endproclaim

As a consequence of the above we obtain the following lemma.
\proclaim{Lemma 2}
 Let $x, a_i, b_{ij}$ all belong to a Banach space $ (B, 
\|\cdot\|)$,
with $b_{ii} =0$. Let $\{\epsilon_i\}$ be a sequence of 
independent and 
symmetric
Bernoulli
random variables, that is, $P(\epsilon_i = 1) = 
P(\epsilon_i =-1) 
={1\over 2}$.
Then,
for a universal constant $c > 0 $,
$$P\left(\left\|x+\sum_{i=1}^n a_i \epsilon_i + 
\sum_{1\le i\ne j \le n} b_{ij} \epsilon_i \epsilon_j 
\right\|\ge \|x\|\right)
\ge c^{-1}.$$
\endproclaim
\demo{Proof} Suppose that $a_i, b_{ij}$ are in $R$; then 
it follows
easily from (1.4) of [9] (see also Sections 6.2 and
6.5 of [10]) that
$$ \left(E\left|\sum_{i=1}^n a_i\epsilon_i + 
   \sum_{1\le i \ne j \le n} b_{ij} \epsilon_i \epsilon_j 
\right|^4\right)^{1\over 4} 
   \le c \left(
   E\left|\sum_{i=1}^n a_i\epsilon_i + 
   \sum_{1\le i \ne j \le n} b_{ij} \epsilon_i \epsilon_j 
\right|^2\right)^{1/2} ,$$ 
for some constant $c>0$.
Next, observe that $\|\xi \|_4 \le c \|\xi \|_2 $ 
implies that 
$\|\xi \|_2 \le c^2 \|\xi \|_1 $ (since $E(\xi )^2 \le 
(E|\xi |)^{2/3}\cdot (E(\xi )^4)^{1/3}$) . 
The result then follows by
Proposition~1.
\enddemo
\demo{Proof 
of Theorem \rm 1}
We first transform the problem of proving (6) into a 
problem dealing  
(conditionally) with a non-homogeneous binomial 
in Bernoulli random variables.
Let 
$\{\epsilon_i \}$ be a sequence of independent
and symmetric
Bernoulli random variables independent of $\{X_i\}$, 
$\{\tilde X_i\}$.  
Let $(Z_i, \tilde Z_i) = (X_i, \tilde X_i) $ if 
$\epsilon_i =  1$ 
and $(Z_i, \tilde Z_i) = (\tilde X_i,  X_i) $ if 
$\epsilon_i = -  1$.
Then,
$$\eqalign{ 4f_{ij}(\tilde Z_i, Z_j) & = 
\{(1-\epsilon_i)(1+\epsilon_j)
f_{ij}(X_i,X_j) 
+ (1+\epsilon_i)(1+\epsilon_j)
f_{ij}(\tilde X_i,X_j) \cr
&\ \ 
+ (1-\epsilon_i)(1-\epsilon_j)
f_{ij}( X_i,\tilde X_j) 
+ (1+\epsilon_i)(1-\epsilon_j)
f_{ij}( \tilde X_i,\tilde X_j) \} .\cr }\leqno (8)$$
Setting ${\Cal G} = \sigma (X_i, \tilde X_i ; i=1,...,n)$,
we get 
$$4E(f_{ij}(\tilde Z_i, Z_j)| {
\Cal G})  =\{ 
f_{ij}(X_i,X_j) + 
f_{ij}(\tilde X_i,X_j) +
f_{ij}( X_i,\tilde X_j) +
f_{ij}( \tilde X_i,\tilde X_j)\}.\leqno (9) $$
From Lemma 2, (3), (8), and (9), and letting $x=T_n$,
it follows that for some $c>0$,
$$P\left(4\left\|\sum_{1\le i\ne j\le n} 
f_{ij}\left(\tilde Z_i, Z_j\right) 
\right\| \ge \|T_n\|
\big| {\Cal G} \right) \ge c^{-1} .
$$
Integrating over the set
$\{\| T_n
\|  \ge t\}$,
we get
$$\split
 {1 \over c}P(\|
T_n\|  \ge t) \le&
P\left(4\left\|\sum_{1\le i\ne j\le n} f_{ij}(\tilde Z_i,  
Z_j) 
\right\|  \ge t\right)\\
 =  &
P\left(4\left\|\sum_{1\le i\ne j\le n} f_{ij}(\tilde X_i,  
X_j) \right\| 
\ge  t\right),  
\endsplit$$ 
%\noindent 
since the sequence\ $\{(X_i,\tilde X_i), i=1,...,n\}$
has the same distribution as
$\{(Z_i,
\widetilde Z_i), i=1,...,n\}.$ The proof is completed by 
using this
inequality along with (4) and (5).

The proof of (2) is similar and  uses an analogue of (8) 
concerning
 $4f(Z_i,Z_j)$. 
In obtaining this bound, one does not need to use Lemma 1. 
Instead
one 
uses the symmetry condition on the functions $f_{ij}$, 
introduced after (1) 
and equation (3), to get
$$\split
&P\left(\left\|\sum_{1\le i \ne j \le n} f_{ij}(\tilde 
X_i,X_j) \right\| \ge t
\right) =
P\left(\left\|\sum_{1\le i \ne j \le n} f_{ij}(X_i,\tilde 
X_j) +
f_{ij}(\tilde X_i,X_j) \right\| \ge 2t\right)  \\
&\qquad \le P(\| T_n
  \| \ge {2\over 3}t)  
+ 2P\left(\left\|\sum_{1\le i \ne j \le n} f_{ij}(X_i,X_j) 
\right\| 
\ge {2\over 3}t\right). \endsplit$$
\enddemo

\enddocument